\documentclass[11pt]{article}
\renewcommand{\baselinestretch}{1.2}
\newcommand{\dated}{\mbox{} \hfill {\small [{\tt \today}]}} \usepackage{amsmath,amssymb,amscd}
\newcommand{\pf}[1]{\trivlist \item[\hskip\labelsep\it #1\ ]}
\newcommand{\varpf}[1]{\trivlist \item[\hskip\labelsep\sc #1:]}
\newcommand{\qedbox}{$\rlap{$\sqcap$}\sqcup$}
\newcommand{\qed}{\qquad \qedbox \endtrivlist}
\newcommand{\varqed}{\hfill \rule{0.6em}{0.6em} \endtrivlist}
\newenvironment{proof}{\pf{Proof}}{\qed}

\newenvironment{remark}{\pf{Remark}}{\endtrivlist}
\newenvironment{remarks}{\pf{Remarks} 
   \begin{enumerate}}{\end{enumerate} \endtrivlist}
\newenvironment{example}{\pf{Example}}{\endtrivlist}
\newenvironment{examples}{\pf{Examples} 
   \begin{enumerate}}{\end{enumerate} \endtrivlist}
\newenvironment{items}{
  \begin{enumerate} 
                    
  }{\end{enumerate}}

\newenvironment{keywords}{\noindent\small {\it Keywords\/}:}{\vskip 4pt}
\newenvironment{classification}{\noindent\small 2000 {\it Mathematics Subject
Classification\/}:}{\vskip 12pt}

\newcommand{\comps}{{\mathbb C}}

\newcommand{\posints}{{\mathbb N}}

\newcommand{\tensor}{\otimes}

\newcommand{\Tensor}{\hat{\otimes}}
\newcommand{\OTensor}{\Hat{\Hat{\otimes}}}

\newcommand{\wTensor}{\check{\otimes}}
\newcommand{\WOTensor}{\Check{\Check{\otimes}}}

\newcommand{\cstar}{{C^\ast}}

\newcommand{\id}{{\mathrm{id}}}

\newcommand{\A}{{\mathfrak A}}

\newcommand{\VN}{\operatorname{VN}}

\newtheorem{theorem}{Theorem}[section]
\newtheorem{lemma}[theorem]{Lemma}
\newtheorem{corollary}[theorem]{Corollary}
\newtheorem{proposition}[theorem]{Proposition}
\newtheorem{df}[theorem]{Definition}
\newenvironment{definition}{\begin{df} \rm}{\end{df}}

\newcommand{\PM}{\operatorname{PM}}
\title{Operator Fig\`a-Talamanca--Herz algebras}
\author{\it Volker Runde\thanks{Research supported by NSERC under grant no.\ 227043-00.}}
\date{}
\begin{document}
\maketitle
\begin{abstract}
Let $G$ be a locally compact group. We use the canonical operator space structure on the spaces $L^p(G)$ for $p \in [1,\infty]$ introduced by G.\ Pisier to define
operator space analogues $OA_p(G)$ of the classical Fig\`a-Talamanca--Herz algebras $A_p(G)$. If $p \in (1,\infty)$ is arbitrary, then $A_p(G) \subset OA_p(G)$ such that
the inclusion is a contraction; if $p = 2$, then $OA_2(G) \cong A(G)$ as Banach spaces spaces, but not necessarily as operator spaces.
We show that $OA_p(G)$ is a completely contractive Banach algebra for each $p \in (1,\infty)$, and that $OA_q(G) \subset OA_p(G)$ completely contractively for amenable $G$ 
if $1 < p \leq q \leq 2$ or $2 \leq q \leq p < \infty$. Finally, we characterize the amenability of $G$ through the existence of a bounded approximate identity in $OA_p(G)$ for
one (or equivalently for all) $p \in (1,\infty)$.
\end{abstract}
\section*{Introduction}
The Fourier algebra $A(G)$ of a locally compact group $G$ was introduced by P.\ Eymard in \cite{Eym0}. If $G$ is abelian with dual group $\Gamma$, then the Fourier transform induces an isometric
isomorphism of $A(G)$ and $L^1(\Gamma)$. Although the Fourier algebra is an invariant for $G$ --- like $L^1(G)$ ---, its Banach algebraic amenability does not correspond well to the amenability of
$G$ --- very much unlike $L^1(G)$: The group $G$ is amenable if and only if $L^1(G)$ is amenable as a Banach algebra (\cite{Joh1}), but there are compact groups, among them $SO(3)$, for which $A(G)$ 
fails to be amenable (\cite{Joh2}; for more on the amenability of $A(G)$, see \cite{For} and \cite{LLW}).
\par
Since $A(G)$ is the predual of the group von Neumann algebra $\VN(G)$, it is an operator space in a natural manner. In \cite{Rua}, Z.-J.\ Ruan introduced a variant of amenability for Banach algebras that 
takes operator spaces structures into account, called {\it operator amenability\/}. He showed that a locally compact group $G$ is amenable if and only if $A(G)$ is operator amenable. Further results
by O.\ Yu.\ Aristov (\cite{Ari}), B.\ E.\ Forrest and P.\ J.\ Wood (\cite{FW}, \cite{Woo}), as well as N.\ Spronk (\cite{Spr}) lend additional support to the belief that homological properties of $A(G)$, such as 
biprojectivity or weak amenability, reflect properties of $G$ much more naturally if the operator space structure is taken into account.
\par
In \cite{Her1}, C.\ Herz introduced, for a locally compact group $G$ and $p \in (1,\infty)$, an $L^p$-analogue of the Fourier algebra, denoted by $A_p(G)$. These algebras are called
{\it Fig\`a-Talamanca--Herz algebras\/}. Since arbitrary Fig\`a-Talamanca--Herz algebras are not preduals of von Neumann algebras, there is --- at the first glance ---
no natural operator space structure for $A_p(G)$, where $G$ is a locally compact group and $p \in (1,\infty) \setminus \{ 2 \}$. 
\par
In \cite{Pis2}, G.\ Pisier used interpolation techniques to equip the spaces $L^p(G)$ for $p \in (1,\infty)$ with a canonical operator space structure. We use this particular operator space structure on the 
$L^p$-spaces to define operator analogues of $OA_p(G)$ of the classical Fig\`a-Talamanca--Herz algebras. If $p = 2$, the Banach spaces $OA_2(G)$ and $A(G)$ are isometrically ismorphic, but for
$G$ compact and non-abelian, $OA_2(G)$ and $A(G)$ are not isometrically isomorphic as operator spaces. For $p \in (1,\infty) \setminus \{ 2 \}$, we only have a contractive inclusion $A_p(G) \subset OA_p(G)$. 
Nevertheless, many of the classical results on Fig\`a-Talamanca--Herz algebras have analogues in the operator space context. We prove that each $OA_p(G)$ is a completely contractive Banach algebra and obtain  
operator analogues of Herz's classical inclusion result for Fig\`a-Talamanca--Herz algebras on amenable groups and of the Leptin--Herz theorem, which characterizes the amenable locally compact groups through
the existence of bounded approximate identities in their Fig\`a-Talamanca--Herz algebras.
\subsubsection*{Acknowledments}
I am grateful to Anselm Lambert for sharing his expertise on Hilbertian operator spaces with me, to Matthias Neufang for bringing reference \cite{Pis2} to my attention, and to 
Gilles Pisier for pointing out an error in an earlier, stronger version of Proposition \ref{cbprop}. Thanks are also due to the referee whose suggestions helped to shorten and to clarify Section 3.
\section{Preliminaries}
\subsection{Fig\`a-Talamanca--Herz algebras}
Let $G$ be a locally compact group, For any function $\phi \!: G \to \comps$, we define $\check{\phi} \!: G \to \comps$ by letting $\check{\phi}(x) := \phi(x^{-1})$ for $x
\in G$. Let $p,q \in (1,\infty)$ be such that $\frac{1}{p} + \frac{1}{q} = 1$. The {\it Fig\`a-Talamanca--Herz algebra\/} $A_p(G)$ consists of those functions $\phi \!: G \to \comps$
such that there are sequences $( \xi_n )_{n=1}^\infty$ in $L^p(G)$ and $( \eta_n )_{n=1}^\infty$ in $L^q(G)$ such that
\begin{equation} \label{Apeq1}
  \sum_{n=1}^\infty \| \xi_n \|_{L^p(G)} \| \eta_n \|_{L^q(G)} < \infty
\end{equation}
and
\begin{equation} \label{Apeq2}
  \phi = \sum_{n=1}^\infty \xi_n \ast \check{\eta}_n.
\end{equation}
The norm $\| \phi \|_{A_p(G)}$ is defined as the infimum over all sums (\ref{Apeq1}) such that (\ref{Apeq2}) holds. It is clear that $A_p(G)$ --- as a quotient of $L^p(G) \Tensor L^q(G)$ ---
is a Banach space. It was shown by C.\ Herz (\cite{Her1}) that $A_p(G)$ is, in fact, a Banach algebra. The case where $p = q = 2$ was previously studied by P.\ Eymard (\cite{Eym0});
in this case $A(G) := A_2(G)$ is called the {\it Fourier algebra\/} of $G$.
\par
For $p \in [1,\infty]$, let $\lambda_p \!: G \to {\cal B}(L^p(G))$ be the regular left representation of $G$ on $L^p(G)$. It is well known that 
$\lambda_p$ extends --- through integration --- to a representation of $M(G)$ --- and thus of $L^1(G)$ --- on $L^p(G)$. Let $p \in (1,\infty)$. The algebra of
{\it $p$-pseudomeasures\/} $\PM_p(G)$ is defined as the $w^\ast$-closure of $\lambda_p(L^1(G))$; it is easy to see that $\lambda_p(M(G)) \subset \PM_p(G)$. Let $q \in (1,\infty)$ be such
that $\frac{1}{p} + \frac{1}{q} = 1$. Then $\PM_p(G) \cong A_q(G)^\ast$ via
\[
  \langle \xi \ast \check{\eta}, T \rangle := \langle T\eta, \xi \rangle \qquad (\xi \in L^q(G), \, \eta \in L^p(G), \, T \in \PM_p(G)). 
\]
If $p =2$, then $\VN(G) := \PM_2(G)$ is a von Neumann algebra, the {\it group von Neumann algebra\/} of $G$. For more information, see \cite{Eym0}, \cite{Eym}, \cite{Her1}, \cite{Her2}, and \cite{Pie}.
\par
The following result was proved by C.\ Herz (\cite[Theorem C]{Her1}):
\begin{proposition} \label{incl}
Let $G$ be an amenable, locally compact group, and let $p,q \in (1,\infty)$ be such that $1 < p \leq q \leq 2$ or $2 \leq q \leq p < \infty$. Then $A_q(G) \subset A_p(G)$ such that the inclusion 
is a contraction with dense range.
\end{proposition}
\par
Related results can be found in \cite{HR} and \cite{Fur}.
\subsection{Operator spaces}
There are now comprehensive sources on operator space theory available
(\cite{ER}, \cite{Pis3}, \cite{Wit}). We shall thus content ourselves with an outline of the basic concepts and results. In our notation, we mostly follow \cite{ER},
except that we use the symbols $\Tensor$ and $\wTensor$ for the projective and injective {\it Banach\/} space tensor product, respectively.
\par
Let $n \in \posints$, and let $E$ be a vector space. We denote the vector space of $n \times n$-matrices with entries from $E$ by ${\mathbb M}_n(E)$. If $E = \comps$,
we simply let ${\mathbb M}_n := {\mathbb M}_n(\comps)$. We always suppose that ${\mathbb M}_n$ is equipped with the operator norm  $| \cdot |_n$ from its canonical action on $n$-dimensional
Hilbert space. Via matrix multiplication, ${\mathbb M}_n$ acts on ${\mathbb M}_n(E)$.
\begin{definition} 
Let $E$ be a vector space. A {\it matricial norm\/} on $E$ is a family $( \| \cdot \|_n )_{n=1}^\infty$ such that $\| \cdot \|_n$
is a norm on ${\mathbb M}_n(E)$ for each $n \in \posints$, such that the following two axioms are satisfied:
\[
  \| \lambda \cdot x \cdot \mu \|_m \leq | \lambda |_n \| x \|_n | \mu |_n \qquad (\lambda, \mu  \in {\mathbb M}_n, \, x \in {\mathbb M}_n(E)),
\]
and
\[
  \left\| \left[ \begin{array}{c|c} x & 0 \\ \hline 0 & y \end{array} \right] \right\|_{n+m} = \max \{ \| x \|_n, \| y \|_m \} \qquad (x \in  {\mathbb M}_n(E), \, y \in {\mathbb M}_m(E)).
\]
\end{definition}
\begin{definition}
A vector space $E$ equipped with a matricial norm $( \| \cdot \|_n )_{n=1}^\infty$ is called a {\it matricially normed space\/}. 
If each space $({\mathbb M}_n(E), \| \cdot \|_n)$ is a Banach space, $E$ is called an (abstract) {\it operator space\/}.
\end{definition}
\begin{remarks}
\item It is easy to see that a matricially normed space such that $(E, \| \cdot \|_1)$ is complete is already an operator space.
\item In our choice of terminology --- reserving the term ``operator space'' for complete spaces ---, we follow \cite{Wit} rather than \cite{ER}. 
\end{remarks}
\begin{example}
Let $\A$ be a $\cstar$-algebra, and let $\| \cdot \|_n$ be the unique $\cstar$-norm on ${\mathbb M}_n(\A)$. Let $E$ be a closed subspace of $\A$. Then $( \| \cdot \|_n |_{{\mathbb M}_n(E)} )_{n=1}^\infty$
is a matricial norm on $E$. Operator spaces of this form are called {\it concrete operator spaces\/}.
\end{example}
\par
Let $E$ and $F$ be matricially normed spaces, let $T \in {\cal B}(E,F)$, and let $n \in \posints$. Then $T^{(n)} \in {\cal B}({\mathbb M}_n(E),{\mathbb M}_n(F))$, the {\it $n$-th amplification\/} of $T$,
is defined as $\id_{{\mathbb M}_n} \tensor T$ with the usual identifications. 
\begin{definition}
Let $E$ and $F$ be a matricially normed spaces. A map $T \in {\cal B}(E,F)$ is called {\it completely bounded\/} if
\[
  \| T \|_{\mathrm{cb}} := \sup_{n \in \posints} \| T^{(n)} \| < \infty.
\]
If $\| T \|_{\mathrm{cb}} \leq 1$, we say that $T$ is a {\it complete contraction\/}, and if $T^{(n)}$ is an isometry for each $n \in \posints$, we call $T$ a {\it complete isometry\/}.
\end{definition}
\par
For any two matricially normed spaces $E$ and $F$, the collection of all completely bounded maps from $E$ to $F$ is denoted by ${\cal CB}(E,F)$. It is routinely checked that 
$({\cal CB}(E,F), \| \cdot \|_{\mathrm{cb}})$ is a normed space which is complete if $F$ is an operator space.
\par
The following fundamental theorem is due to Z.-J.\ Ruan (\cite{Rua1}; \cite[Theorem 2.3.5]{ER}):
\begin{theorem}
Let $E$ be an operator space. Then $E$ is completely isometrically isomorphic to a concrete operator space.
\end{theorem}
\par
Let $E$ and $F$ be matrically normed spaces. Then,  through the canonical identifications
\[
  {\mathbb M}_n({\cal CB}(E,F)) = {\cal CB}(E,{\mathbb M}_n(F)) \qquad(n \in \posints),
\]
the space ${\cal CB}(E,F)$ becomes again a matricially normed space. If $F = \comps$, then $E^\ast$ is isometrically isomorphic to ${\cal CB}(E,\comps)$ (\cite[Proposition 2.2.2]{ER}), and
thus is an operator space in a canonical manner. In particular, dual spaces of $\cstar$-algebras and, more generally, predual spaces of von Neumann algebras can be equipped with an operator space structure
in this particular way.
\par
We shall also extensively use the interpolation techniques for operator spaces developed by G.\ Pisier and expounded in \cite{Pis1} and \cite{Pis2}.
\par
Let $(E_0,E_1)$ be a compatible couple of Banach spaces in the sense of interpolation theory (\cite[2.3]{BL}). For $\theta \in (0,1)$, let $E_\theta := (E_0,E_1)_\theta$ denote the space
obtained from A.\ P.\ Calder\'on's complex method of interpolation (\cite[Chapter 4]{BL}). For each $n \in \posints$, one defines a Banach space structure on ${\mathbb M}_n(E_\theta)$
by letting 
\[
  {\mathbb M}_n(E_\theta) := ({\mathbb M}_n(E_0), {\mathbb M}_n(E_1))_\theta.
\]
In this way, $E_\theta$ becomes an operator space (see \cite[\S2]{Pis1}).
\par
For more details, see \cite{Pis1} and \cite{Pis2}.
\section{Operator space analogues of Fig\`a-Talamanca--Herz algebras}
If $G$ is a locally compact group, $\VN(G) \subset {\cal B}(L^2(G))$ is a concrete operator space. As the predual of $\VN(G)$, the Fourier algebra carries a natural 
operator space structure. For $A_p(G)$ with $p \in (1,\infty) \setminus \{ 2 \}$, no such canonical operator space structure is so obviously available. In this section, we identify
a natural matricial norm on the algebra $\lambda_p(M(G))$ of convolution operators on $L^p(G)$, which, in turn, enables us to define an operator analogue of $A_p(G)$.
\par
Let $G$ be a locally compact group. As a commutative von Neumann algebra $L^\infty(G)$ has a canonical operator space structure, and so does its 
predual $L^1(G)$.
The couple $(L^\infty(G), L^1(G))$ of Banach spaces is compatible in the sense of interpolation theory (see \cite{BL}). For each $\theta \in (0,1)$, the 
complex interpolation space $(L^\infty(G),L^1(G))_\theta$ is thus well defined; for $p \in (1,\infty)$ and $\theta = \frac{1}{p}$, it is well known that
\begin{equation} \label{isom}
  (L^\infty(G),L^1(G))_\theta \cong L^p(G)
\end{equation}
isometrically (\cite[Theorem 5.1.1]{BL}). In \cite{Pis1}, G.\ Pisier shows that complex interpolation between Banach spaces carrying an operator space structure can be used to equip
the resulting interpolation spaces with operator space structures. In view of (\ref{isom}), this can be used to define an operator space structure on each of the spaces $L^p(G)$ with 
$p \in [1, \infty]$ (\cite{Pis2}); we denote this operator space by $OL^p(G)$.
Since $L^\infty(G)$ is a commutative $\cstar$-algebra, $OL^\infty(G) = \min L^\infty(G)$ holds (\cite[Proposition 3.3.1]{ER}). In particular, 
each $T \in {\cal B}(L^\infty(G))$ is completely bounded such that
\[
  \| T \|_{\mathrm{cb}} = \| T\| \qquad (T \in {\cal B}(L^\infty(G))).
\]
In the terminology of \cite{Pis1}, $OL^\infty(G)$ is homogeneous. Since $OL^1(G)^\ast = OL^\infty(G)$, the same is true for $OL^1(G)$. 
Let $\mu \in M(G)$, and let $p \in (1,\infty)$. Since $\lambda_p(\mu)$ is obtained from $\lambda_\infty(\mu)$ and $\lambda_1(\mu)$ through 
interpolation, it follows from \cite[Proposition 2.1]{Pis1} that $\lambda_p(\mu) \in {\cal CB}(OL^p(G))$ such that $\| \lambda_p(\mu) \|_{\mathrm{cb}} \leq \| \mu \|$. 
We can say even more for particular $p$ and $\mu$:
\begin{proposition} \label{cbprop}
Let $G$ be a locally compact group, let $p \in (1,\infty)$, and let $\mu \in M(G)$. Then $\lambda_p(\mu) \in {\cal CB}(OL^p(G))$ such that
\[
  \| \lambda_p(\mu) \|_{\mathrm{cb}} \leq \| \mu \|.
\]
If $p =2$ or if $\mu$ is positive, we even have
\begin{equation} \label{cbeq}
  \| \lambda_p(\mu) \|_{\mathrm{cb}} = \| \lambda_p(\mu) \|.
\end{equation}
\end{proposition}
\begin{proof}
The first part is clear in view of the remarks made immediately before.
\par
If $p =2$, then $OL^2(G)$ is completely isometrically isomorphic to the operator Hilbert space $OH(I)$ for an appropriate index set $I$ (\cite[Proposition 2.1(iii)]{Pis2}).
Since $OH(I)$ is homogeneous by \cite[Proposition 1.5(i)]{Pis1}, this establishes (\ref{cbeq}) for $p =2$.
\par
Let $p \in (1,\infty)$ be arbitrary, and let ${\cal S}^p$ denote the $p$-th von Neumann--Schatten class on $\ell^2$. By \cite[Proposition 2.4]{Pis2}, $\| \lambda_p(\mu) \|_{\mathrm{cb}}$ equals
the operator norm of $\lambda_p(\mu) \tensor \id_{{\cal S}^p}$ on $L^p(G,{\cal S}^p)$. If $\mu$ is positive, then $\lambda_p(\mu)$ is a positive operator on $L^p(G)$. Hence,
$\| \lambda_p(\mu) \tensor \id_{{\cal S}^p} \| = \| \lambda_p(\mu) \|$ holds by \cite[7.3, Theorem]{DF}.
\end{proof}
\begin{remarks}
\item In an earlier preprint version of this paper, we claimed that $\| \lambda_p(\mu) \|_{\mathrm{cb}} = \| \lambda_p(\mu) \|$ for all $p \in (1,\infty)$ and for all $\mu \in M(G)$. This is false as G.\ Pisier 
pointed out to us. In fact, for $p \neq 2$, the norms of ${\cal B}(L^p(G))$ and ${\cal CB}(OL^p(G))$ need not even be equivalent on $\lambda_p(M(G))$. {\it If\/} they were equivalent, then 
\cite[Th\'eor\`eme 2.2]{Eym} could be used to show that $\PM_p(G) \subset {\cal CB}(OL^p(G))$ (with equivalent norms) for amenable $G$. For $G$ abelian, taking the Fourier transform would thus yield that 
every Fourier multiplier on $L^p(G)$ is completely bounded. For compact, infinite, abelian $G$, however, there are counterexamples (\cite[Proposition 8.1.3]{Pis2}).
\item Using (\ref{cbeq}) for $p =2$ and interpolating between $OL^\infty(G)$ and $OL^2(G)$ (or $OL^2(G)$ and $OL^1(G)$, depending on whether $p \geq 2$ or $p \leq 2$), we obtain slightly better
estimates for $\| \lambda_p(\mu) \|_{\mathrm{cb}}$ than $\| \mu \|$.
\end{remarks}
\par
The following lemma is a simple, special case of \cite[Theorem 4.1]{Pis2}:
\begin{lemma}
Let $G$ be a locally compact group, and let $p,q \in (1,\infty)$ be such that $\frac{1}{p} + \frac{1}{q} = 1$. Then $OL^p(G)^\ast \cong OL^q(G)$ as operator spaces.
\end{lemma}
\par
We denote the projective tensor product of operator spaces by $\OTensor$. Let $p,q \in (1,\infty)$ be such that $\frac{1}{p} + \frac{1}{q} = 1$.
By \cite[Proposition 7.1.4 and Corollary 7.1.15]{ER}, we have a completely isometric isomorphism $(OL^p(G) \OTensor OL^q(G))^\ast \cong {\cal CB}(OL^q(G))$. Hence, 
the following definition is meaningful:
\begin{definition} \label{FTHdef}
Let $G$ be a locally compact group, and let $p,q \in (1,\infty)$ be such that $\frac{1}{p} + \frac{1}{q} = 1$. Then the {\it operator Fig\`a-Talamanca--Herz\/} algebra
$OA_p(G)$ is defined as the quotient of $OL^p(G) \OTensor OL^q(G)$ under the restriction to $\lambda_q(M(G))$.
\end{definition}
\begin{remarks}
\item It makes no difference if, in Definition \ref{FTHdef}, we replace $\lambda_q(M(G))$ by the smaller space $\lambda_q(L^1(G))$: Since the closed unit ball of $L^1(G)$ is $w^\ast$-dense in 
the closed unit ball of $M(G)$, the space $\lambda_q(L^1(G))$ is $w^\ast$-dense in $\lambda_q(M(G))$.
\item The quotient map from $OL^p(G) \OTensor OL^q(G)$ onto $OA_p(G)$ is the linearization of the bilinear map $L^p(G) \times L^q(G) \ni (f,g) \mapsto f \ast \check{g}$.
\item Since $\lambda_q \!: M(G) \to {\cal CB}(OL^q(G))$ is a $w^\ast$-continuous complete contraction, taking adjoints yields an injective complete contraction from $OA_p(G)$ into ${\cal C}_0(G)$.
\item Since $L^p(G) \Tensor L^q(G)$ embeds contractively into $OL^p(G) \OTensor OL^q(G)$, we have a canonical, necessarily contractive inclusion $A_p(G) \subset OA_p(G)$.
\end{remarks}
\par
The following is the operator space analogue of a classical result on Fig\`a-Talamanca--Herz algebras:
\begin{proposition} \label{checkprop}
Let $G$ be a locally compact group, and let $p, q \in (1,\infty)$ be such that $\frac{1}{p} + \frac{1}{q} = 1$. Then
\begin{equation} \label{checkmap}
  OA_p(G) \to OA_q(G), \quad \phi \mapsto \check{\phi}
\end{equation}
is a completely isometric isomorphism.
\end{proposition}
\begin{proof}
By \cite[Proposition 7.1.4]{ER}, the flip map
\begin{equation} \label{flip}
  OL^p(G) \OTensor OL^q(G) \to OL^q(G) \OTensor OL^p(G), \quad \xi \tensor \eta \mapsto \eta \tensor \xi
\end{equation}
is a completely isometric isomorphism. Since the diagram
\[
  \begin{CD}  
  OL^p(G) \OTensor OL^q(G) @>(\ref{flip})>>     OL^q(G) \OTensor OL^p(G) \\
  @VVV                                          @VVV                     \\
  OA_p(G)                  @>(\ref{checkmap})>> OA_q(G)
  \end{CD}
\]
commutes, this and Definition \ref{FTHdef} establish the claim.
\end{proof}
\section{Operator $p$-spaces}
In this section, we leave the framework of locally compact groups, and work with $L^p$-spaces over arbitrary measure spaces $X$. Although $L^\infty(X)$ then need no longer be
a von Neumann algebra, it is still a commutative $\cstar$-algebra, and $L^1(X)$ embeds isometrically into $L^\infty(X)^\ast$. We can thus still meaningfully speak of the operator
$L^p$-spaces $OL^p(X)$ for $p \in [1,\infty]$ (\cite{Pis2}).
\par
Let $X$ be a measure space, and let $E$ be an arbitrary operator space. Define
\[
  OL^1(X,E) := OL^1(X) \OTensor E.
\]
Suppose that $E$ is represented as a closed subspace of ${\cal B}({\mathfrak H})$ for some Hilbert space. Then $L^\infty(X,E)$ is a closed subspace of the $\cstar$-algebra
$L^\infty(X,{\cal B}({\mathfrak H}))$ and thus an operator space --- denoted by $OL^\infty(X,E)$ --- in a canonical manner. Interpolating and letting
\[
  OL^p(X,E) := (OL^\infty(X,E), OL^1(X,E))_\theta
\]
for $p \in (1,\infty)$ and $\theta = \frac{1}{p}$, we equip the spaces $L^p(X,E)$ with operator space structures (see \cite[5.1.2 Theorem]{BL}). 
\par
Let $\WOTensor$ denote the injective tensor product of operator spaces, and define
\[
  OL^\infty_0(X,E) := OL^\infty(X) \WOTensor E.
\]
Then $OL^\infty_0(X,E)$ can be identified with a closed subspace of $OL^\infty(X,E)$. In general, we have $OL^\infty_0(X,E) \subsetneq OL^\infty(X,E)$, but an exhaustion argument and \cite[Lemma 0.1]{Pis2} applied to each matrix level yield:
\begin{lemma} \label{exhaust}
Let $X$ be a measure space, let $E$ be an operator space, and let $p \in (1,\infty)$. Then, with $\theta = \frac{1}{p}$, we have a completely isometric isomorphism
\[
  OL^p(X,E) = (OL^\infty_0(X,E), OL^1(X,E))_\theta.
\]
\end{lemma}
\begin{remark}
Alternatively, Lemma \ref{exhaust} can be deduced from \cite[Theorem 4.2.2]{BL} (applied to each matrix level).
\end{remark}
\par
Let $E$ and $F$ be matricially normed spaces. A matricial norm on $E \tensor F$ is called a {\it matricial subcross norm\/} (\cite[p.\ 124]{ER}) if
\begin{equation} \label{subcrossdef}
  \| x \tensor y \|_{n_1 n_2} \leq \| x \|_{n_1} \| y \|_{n_2} \qquad (n_1, n_2 \in \posints, \, x \in {\mathbb M}_{n_1}(E), \, y \in {\mathbb M}_{n_2}(F)).
\end{equation}
If equality holds in (\ref{subcrossdef}), we speak of a {\it matricial cross norm\/}.
\begin{proposition} \label{subprop}
Let $X$ be a measure space, let $p \in [1,\infty]$, and let $E$ be any operator space. Then the restriction of the matricial norm on $OL^p(X,E)$ to $L^p(X) \tensor E$ is a matricial cross
norm, and there is a canonical complete contraction from $OL^p(X) \OTensor E$ to $OL^p(X,E)$.
\end{proposition}
\begin{proof}
The claim is true for $p \in \{1, \infty\}$ because the projective and the injective operator tensor norms are matricial cross norms. For arbitrary $p \in [1,\infty]$, it follows through interpolation
that we are dealing with matricial subcross norms on $L^p(X) \tensor E$.
\par
Let $p,q \in (1,\infty)$ be such that $\frac{1}{p} + \frac{1}{q} = 1$. Since $OL^q(X,E^\ast)$ embeds completely isometrically into $OL^p(X,E)^\ast$ by \cite[Theorem 4.1]{Pis2},
and since the matricial norm on $L^q(X) \tensor E^\ast$ is matricially subcross, an elementary calculation shows that the matricial norm on $L^p(X) \tensor E$ is, in fact,
a matricial cross norm.
\par
The claim about a canonical complete contraction from $OL^p(X) \OTensor E$ to $OL^p(X,E)$ is then an immediate consequence of \cite[Theorem 7.1.1]{ER}.
\end{proof}
\par
Our next proposition is well known in the Banach space setting (see \cite[7.3]{DF}):
\begin{proposition} \label{mapprop2}
Let $X$ be a measure space, let $p \in [1,\infty)$, and let $E$ be any operator space. Then the amplification map
\begin{equation} \label{ampl0}
  {\cal CB}(E) \ni T \mapsto \id_{L^p(X)} \tensor T
\end{equation}
is a complete isometry from ${\cal CB}(E)$ to ${\cal CB}(OL^p(X,E))$.
\end{proposition}
\begin{proof}
As a consequence of Proposition \ref{subprop}, it is sufficient to prove that (\ref{ampl0}) is a complete contraction.
Since both the projective and the injective matricial norm are uniform operator space tensor norms (\cite[Proposition 5.11]{BlP}; see also footnote 26 of \cite{Wit}) in the sense of \cite[Definition 5.9]{BlP},
this is clear for $p \in \{ 1, \infty \}$ (in the case where $p =\infty$, replace $OL^\infty(X,E)$ by $OL^\infty_0(X,E)$). 
\par
For $p \in (1,\infty)$, the claim then follows from \cite[Proposition 2.1]{Pis1}.
\end{proof}
\begin{remark}
That (\ref{ampl0}) is a contraction is observed on \cite[p.\ 39]{Pis2}.
\end{remark}
\par
The question of whether the r\^oles of $OL^p(X)$ and $E$ in Proposition \ref{mapprop2} can be interchanged is at the heart of the following definition, which is the operator space analogue 
of the central concept of \cite{Her1}:
\begin{definition} \label{pdef}
Let $p \in (1,\infty)$. An operator space $E$ is called an {\it operator $p$-space\/} if, for each measure space $X$, the amplification map
\begin{equation} \label{ampl}
  {\cal CB}(OL^p(X)) \ni T \mapsto T \tensor \id_E 
\end{equation}
is a complete isometry from ${\cal CB}(OL^p(X))$ to ${\cal CB}(OL^p(X,E))$. 
\end{definition}
\begin{remarks}
\item If $E$ is not an operator $p$-space, then it is not even clear that (\ref{ampl}) attains its values in ${\cal CB}(OL^p(X,E))$.
\item The definition of an operator $p$-space is reminiscent of the equivalent conditions in \cite[Theorem 6.9]{Pis2} (in the case where $p = 2$; see also \cite[Corollary 7.2.7]{Pis2}). We are, 
however, interested in whether or not (\ref{ampl}) is a complete isometry, not just an isometry.
\end{remarks}
\par
From \cite[(2.5)]{Pis2} and \cite[Lemma 1.1]{Ble}, we obtain immediately:
\begin{corollary} \label{pcor}
Let $p,q \in (1,\infty)$ be such that $\frac{1}{p} + \frac{1}{q} = 1$, and let $E$ be an operator $p$-space. Then $E^\ast$ is an operator $q$-space.
\end{corollary}
\par
\begin{theorem} \label{pthm0}
Let $X$ be a measure space, and let $p \in (1,\infty)$. Then $OL^p(X)$ is an operator $p$-space. 
\end{theorem}
\begin{proof}
In view of Proposition \ref{mapprop2} and \cite[(3.6)]{Pis2}, this is obvious.
\end{proof}
\par
The class of operator $p$-spaces is stable under complex interpolation:
\begin{lemma} \label{ipoll}
Let $p \in (1,\infty)$, and let $(E_0,E_1)$ be a compatible couple of operator spaces such that $E_0$ and $E_0$ are operator $p$-spaces. Then $(E_0,E_1)_\theta$ is an operator $p$-space
for each $\theta \in (0,1)$. 
\end{lemma}
\begin{proof}
Let $X$ be a measure space. Since
\[
  OL^p(X,(E_0,E_1)_\theta) \cong (OL^p(X,E_0), OL^p(X,E_1))_\theta \qquad (\theta \in (0,1))
\]
by \cite[(2.1)]{Pis2}, the claim is immediate from \cite[Proposition 2.1]{Pis1}.
\end{proof}
\par
Let $\kappa$ be any cardinal number, and let $p \in [1,\infty]$. We use ${\cal S}^p_\kappa$ to denote the $p$-th von Neumann--Schatten class on a Hilbert space of dimension $\kappa$ (so that
${\cal S}^p_{\aleph_0} = {\cal S}^p$). For $p \in (1,\infty)$ and $\theta = \frac{1}{p}$, we have
\begin{equation} \label{interp}
  {\cal S}^p_\kappa = ({\cal S}^\infty_\kappa, {\cal S}^1_\kappa )_\theta
\end{equation}
Since ${\cal S}^\infty_\kappa$ and ${\cal S}^1_\kappa$ are operator spaces ${\cal OS}^\infty_\kappa$ and ${\cal OS}^1_\kappa$ in a natural 
fashion, we can use (\ref{interp}), to equip ${\cal S}^p_\kappa$ with an operator space structure. More generally, for any operator space $E$, let
\[
  {\cal S}_\kappa^\infty(E) := {\cal S}_\kappa^\infty \WOTensor E \qquad\text{and}\qquad {\cal S}_\kappa^1(E) := {\cal S}_\kappa^1 \OTensor E ,
\]
and define, for $p \in (1,\infty)$ and $\theta = \frac{1}{p}$:
\[
  {\cal OS}^p_k(E) := ({\cal S}^\infty_\kappa(E), {\cal S}^1_\kappa(E) )_\theta.
\]
See \cite{Pis2}, for more information. It is not hard to see that an operator space $E$ is an operator $p$-space if the amplification map
\[
  {\cal CB}({\cal OS}^p_n) \to {\cal CB}({\cal OS}^p_n(E)), \quad T \mapsto T \tensor \id_E
\]
is a complete contraction for each $n \in \posints$. Consequently, by \cite[(3.6)]{Pis2}, ${\cal S}^p_\kappa$ is an operator $p$-space for each $\kappa$.
\par
For any Hilbert space $H$ --- with Hilbert space dimension $\kappa$ --- let $OH$ denote the corresponding operator Hilbert space as introduced by G.\ Pisier in \cite{Pis1}. Let $R_p$ and $C_p$ denote the rows and colums of ${\cal OS}^p_\kappa$,
respectively. Since $OH = (R_p,C_p)_{\frac{1}{2}}$, and since $R_p$ and $C_p$ --- as completely $1$-complemented subspaces of ${\cal OS}^p_\kappa$ --- are operator $p$-spaces, Lemma \ref{ipoll} yields:
\begin{theorem} \label{pthm}
Let $p \in (1,\infty)$, and let $H$ be a Hilbert space. Then $OH$ is an operator $p$-space.
\end{theorem}
\par
When we return to the setting of locally compact groups in the next section, we require the following corollary of Theorem \ref{pthm}, which is the operator space version of
\cite[Theorem 1]{Her1}:
\begin{corollary} \label{pcor2}
Let $X$ be any measure space, and let $1 < p \leq q \leq 2$ or $2 \leq q \leq p < \infty$. Then $OL^q(X)$ is an operator $p$-space. 
\end{corollary}
\begin{proof}
By Theorem \ref{pthm0}, the claim is clear for $q =p$. For $q =2$, it is an immediate consequence of Theorem \ref{pthm} and \cite[Proposition 2.1(iii)]{Pis2}. For all other $q$ such that
$1 < p \leq q \leq 2$ or $2 \leq q \leq p < \infty$, the claim follows through interpolation from Lemma \ref{ipoll}.
\end{proof}
\section{Operator Fig\`a-Talamanca--Herz algebras as completely contractive Banach algebras}
Let $G$ be a locally compact group. We have {\it called\/} the spaces $OA_p(G)$ for $p \in (1,\infty)$ operator Fig\`a-Talamanca--Herz algebras without bothering 
about whether they are indeed algebras. By Proposition \ref{cbprop} and Kaplansky's density theorem, we have $A(G) \cong OA_2(G)$ as Banach spaces, so that $OA_2(G)$ is a Banach algebra, but it is unclear
how multiplicative and operator space structure interact.
\par
In this section, we shall see that general operator Fig\`a-Talamanca--Herz algebras are not only Banach algebras, but are completely contractive in the following sense
(see \cite{Rua} and \cite{ER}):
\begin{definition}
A Banach algebra $\A$ equipped with an operator space structure is called {\it completely contractive\/} if the algebra product
\begin{equation} \label{mult}
  \A \times \A \to \A, \quad (a,b) \mapsto ab
\end{equation}
is a completely contractive bilinear map. 
\end{definition}
\begin{examples}
\item For any Banach algebra $\A$, its maximal operator space $\max \A$ is a completely contractive Banach algebra.
\item Let $\mathfrak H$ be a Hilbert space. Then every closed subalgebra of ${\cal B}({\mathfrak H})$ is completely contractive.
\item If $G$ is a locally compact group, then $A(G)$ equipped with its canonical operator space structure is a completely contractive
Banach algebra (\cite{Rua}). Note that for infinite, abelian $G$, the Fourier algebra $A(G) \cong L^1(\Gamma)$ is not Arens regular
and thus cannot be isomorphic to a closed subalgebra of ${\cal B}({\mathfrak H})$ for any Hilbert space $\mathfrak H$.
\end{examples}
\par
The following extends Proposition \ref{ccprop}:
\begin{proposition} \label{mixprop}
Let $G$ and $H$ be locally compact groups, and let $p,q \in [1, \infty]$. Then we have $\lambda_p(\mu) \tensor \lambda_q(\nu) \in {\cal CB}(OL^p(G,OL^q(H)))$ for all
$\mu \in M(G)$ and $\nu \in M(H)$ such that 
\[
  \| \lambda_p(\mu) \tensor \lambda_q(\nu) \|_{\mathrm{cb}} \leq \| \mu \| \| \lambda_q(\nu) \|_{\mathrm{cb}}.
\]
\end{proposition}
\begin{proof}
For $p \in \{1,\infty\}$ the claim is an immediate consequence of Proposition \ref{ccprop} and of the mapping properties of the projective and the injective tensor product of operator spaces.
For $p \in (1,\infty)$, it follows from \cite[Proposition 2.1]{Pis2}.
\end{proof}
\begin{remark}
Let $G$ and $H$ be locally compact groups, and let $p,q \in (1,\infty)$ be such that $\frac{1}{p} + \frac{1}{q} = 1$. By Proposition \ref{mixprop}, we have a canonical contractive representation 
$\lambda_{p,q}$ of $L^1(G) \Tensor L^1(H) \cong L^1(G \times H)$ in ${\cal CB}(OL^p(G,OL^q(H)))$. Since $L^1(G \times H)$ is an ideal in $M(G \times H)$ and has an approximate identity bounded by one,
$\lambda_{p,q}$ extends canonically to $M(G \times H)$ as a (completely) contractive representation. 
\end{remark}
\par
Let $G$ and $H$ be locally compact groups, and let $p,q \in (1,\infty)$ be such that $\frac{1}{p} + \frac{1}{q} = 1$. By \cite[Theorem 4.1]{Pis2}, we have
$OL^p(G,OL^q(H))^\ast \cong OL^q(G,OL^p(H))$, so that the following definition is meaningful:
\begin{definition} \label{FThdef2}
Let $G$ and $H$ be locally compact groups, let $p,q \in (1, \infty)$, and let $r,s \in (1, \infty)$ be such that $\frac{1}{p} + \frac{1}{r} = 1$ and 
$\frac{1}{q} + \frac{1}{s} = 1$. Then $OA_{p,q}(G \times H)$ is defined to be the quotient of $OL^p(G,L^q(H)) \OTensor OL^r(G,OL^s(H))$
under the restriction to $\lambda_{r,s}(M(G \times H))$.
\end{definition}
\begin{remarks}
\item By \cite[(3.6)]{Pis2}, it is clear that $OA_{p,p}(G \times G) = OA_p(G \times G)$.
\item We can replace $\lambda_{r,s}(M(G \times H))$ in Definition \ref{FTHdef} by $\lambda_{r,s}(L^1(G \times H))$
\item As for the operator Fig\`a-Talamanca--Herz algebras, it is easy to see that $OA_{p,q}(G \times H)$ embeds completely contractively into ${\cal C}_0(G \times H)$.
\end{remarks}
\begin{proposition} \label{ccprop}
Let $G$  and $H$ be a locally compact groups, and let $p,q \in (1,\infty)$. Then there is a canoncial complete contraction
from $OA_p(G) \OTensor OA_q(H)$ into $OA_{p,q}(G \times H)$.
\end{proposition}
\begin{proof}
Choose $r,s \in (1,\infty)$ such that $\frac{1}{p} + \frac{1}{r} =1$ and $\frac{1}{q} + \frac{1}{s} = 1$. We have a canonical completely isometric isomorphism 
\begin{equation} \label{shuffle}
   (OL^p(G) \OTensor OL^r(G)) \OTensor (OL^q(H) \OTensor OL^s(H)) \cong (OL^p(G) \OTensor OL^q(H)) \OTensor (OL^r(G) \OTensor OL^s(H))
\end{equation}
by \cite[Proposition 7.1.4]{ER}. Consider the diagram
\[
  \begin{CD}
  (OL^p(G) \OTensor OL^r(G)) \OTensor (OL^q(H) \OTensor OL^s(H)) @>>> OL^p(G,OL^q(H)) \OTensor OL^r(G,OL^s(H)) \\
  @VVV                                                                @VVV \\
  OA_p(G) \OTensor OA_p(H)                                       @>>> OA_{p,q}(G \times H),
  \end{CD}
\]
where the top row is the composition of (\ref{shuffle}) with the canonical complete contractions from $OL^p(G) \OTensor OL^q(H)$ to $OL^p(G,OL^q(H))$ and from
$OL^r(G) \OTensor OL^s(H)$ to $OL^r(G,OL^s(H))$, which exist according to Proposition \ref{subprop}. Clearly, going 
along the top row and down the second column is a complete contraction, that factors through the kernel of the first column. Since the first column is a quotient map by 
Definition \ref{FTHdef} and \cite[Proposition 7.1.7]{ER}, we obtain the bottom row, which makes the diagram commutative and clearly is a complete contraction.
\end{proof}
\begin{remark}
In the special case where $G = H$ and $p = q$, we have a canonical complete contraction from $OA_p(G) \OTensor OA_p(G)$ into $OA_p(G \times G)$.
\end{remark}
\par
We require two further lemmas:
\begin{lemma} \label{cclem1}
Let $G$ and $H$ be locally compact groups, and let $1 < p \leq q \leq 2$ or $2 \leq q \leq p < \infty$. Then
\[
  \lambda_p(M(G)) \to {\cal CB}(OL^p(G,OL^q(H))), \quad \lambda_p(\mu) \mapsto \lambda_p(\mu) \tensor \id_{L^q(H)}
\] 
is a complete isometry.
\end{lemma}
\begin{proof}
By Corollary \ref{pcor2}, $OL^q(H)$ is an operator $p$-space. 
\end{proof}
\par
Let $G$ be any locally compact group, and define the {\it fundamental operator\/} $W \!: \comps^{G \times G} \to \comps^{G \times G}$ by letting
\begin{equation} \label{fund1}
  (W \xi)(x,y) := \xi(x,xy) \qquad (\xi \in \comps^{G \times G} , \, x,y \in G).
\end{equation}
Then $W$ is easily seen to be bijective with $W^{-1} \!: \comps^{G \times G} \to \comps^{G \times G}$ being given by
\begin{equation} \label{fund2}
  (W^{-1} \xi)(x,y) := \xi(x,x^{-1}y) \qquad (\xi \in  \comps^{G \times G}, \, x,y \in G).
\end{equation}
\begin{lemma} \label{cclem2}
Let $G$ be a locally compact group, and let $p \in [1,\infty]$. Then $W |_{OL^p(G,OL^q(G))}$ induces a completely isometric isomorphism
$W_{p,q} \!: OL^p(G,OL^q(G)) \to OL^p(G,OL^q(G))$.
\end{lemma}
\begin{proof}
Obviously, $W_{1,1}$, $W^{-1}_{1,1}$, $W_{\infty,\infty}$, and $W_{\infty,\infty}^{-1}$ are complete contractions and thus complete isometries.
\par
Fix $n \in \posints$. Since
\[
  {\mathbb M}_n(OL^\infty_0(G,OL^1(G))) \cong {\mathbb M}_n \WOTensor OL^\infty(G) \WOTensor OL^1(G) \cong OL^\infty(G) \WOTensor {\mathbb M}_n(OL^1(G)),
\]
we have an isometric isomorphism
\[
  {\mathbb M}_n(L^\infty_0(G,L^1(G))) \cong L^\infty(G) \wTensor {\mathbb M}_n(L^1(G)),
\]
which, by $w^\ast$-continuity, extends to an isometric isomorphism
\[
  {\mathbb M}_n(L^\infty(G,L^1(G))) \cong L^\infty(G,{\mathbb M}_n(L^1(G)))
\]
at the Banach space level. A straightforward calculation --- involving left invariance of Haar measure --- then shows that $W_{\infty,1}^{(n)}$ is an isometry. Hence, $W_{\infty,1}$ is a complete isometry.
Analoguously, the claim for $W^{-1}_{\infty,1}$ is established. Since 
\[
  W_{1,\infty} = (W^{-1}_{\infty,1})^\ast |_{OL^1(G,OL^\infty(G))} \qquad\text{and}\qquad W_{1,\infty}^{-1} = (W_{\infty,1})^\ast |_{OL^1(G,OL^\infty(G))}, 
\]
we conclude that $W_{1,\infty}$ and $W_{1,\infty}^{-1}$ are also complete isometries.
\par
It is then an immediate consequence of \cite[Proposition 2.1]{Pis1}, that $W_{p,q}$ and $W^{-1}_{p,q}$ are complete contractions --- and thus complete isometries ---
for all $p \in [1,\infty]$ and $q \in \{1,\infty\}$.
\par
Using \cite[(2.1)]{Pis2} and \cite[Proposition 2.1]{Pis1} again, we obtain the claim for all $p,q \in [1,\infty]$.
\end{proof}
\begin{theorem} \label{ccthm}
Let $G$ be a locally compact group, and let $1 < p \leq q \leq 2$ or $2 \leq q \leq p < \infty$. Then pointwise multiplication induces a complete contraction
from $OA_q(G) \OTensor OA_p(G)$ into $A_p(G)$.
\end{theorem}
\begin{proof}
Let $r,s \in (1,\infty)$ be such that $\frac{1}{p} + \frac{1}{r} = 1$ and $\frac{1}{q} + \frac{1}{s} = 1$. By Lemmas \ref{cclem1} and \ref{cclem2}, the map
\[
  \nabla \!: \lambda_r(M(G)) \to {\cal CB}(OL^r(G,OL^s(G)), \quad \lambda_r(\mu) \mapsto W_{r,s}^{-1}(\lambda_r(\mu) \tensor \id_{L^s(G)})W_{r,s}
\]
is a complete isometry. It is routinely checked that $\nabla(\lambda_r(M(G)) \subset \lambda_{r,s}(M(G \times G))$. It is also immediate that $\nabla$ 
is continuous with respect to the $w^\ast$-topologies on ${\cal CB}(OL^r(G))$ and ${\cal CB}(OL^r(G,OL^s(G)))$, respectively, so that
$\nabla^\ast(OA_{p,q}(G \times G)) \subset OA_p(G)$. Another routine calculation shows that
\[
  (\nabla^\ast \phi)(x) = \phi(x,x) \qquad (\phi \in OA_{p,q}(G \times G), x \in G),
\]
i.e.\ $\nabla^\ast$ restricts functions in $OA_p(G \times G)$ to the diagonal $\{ (x,x) : x \in G \} \subset G \times G$. As the adjoint of a complete isometry,
$\nabla^\ast$ is a complete quotient map and thus, in particular, a complete contraction. Together with Proposition \ref{ccprop}, we thus obtain a completely contractive map
\[
  OA_p(G) \OTensor OA_q(G) \to OA_{p,q}(G \times G) \to OA_p(G), \quad \phi \tensor \psi \mapsto \phi \psi. 
\]
This completes the proof.
\end{proof}
\par
Applying Theorem \ref{ccthm} to the case where $p = q$, we obtain:
\begin{corollary}
Let $G$ be a locally compact group, and let $p \in (1,\infty)$. Then $OA_p(G)$, with pointwise multiplication, is a completely contractive Banach algebra.
\end{corollary}
\par
Another consequence of Theorem \ref{ccthm} is the following operator version of Proposition \ref{incl}:
\begin{corollary} \label{inclcor}
Let $G$ be an amenable, locally compact group, and let $1 < p \leq q \leq 2$ or $2 \leq q \leq p < \infty$. Then $OA_q(G) \subset OA_p(G)$ such that the inclusion is a complete
contraction with dense range.
\end{corollary}
\begin{proof}
By \cite[Theorem 10.4]{Pie}, $A_p(G)$ has an approximate identity bounded by one, say $( e_\alpha )_\alpha$. The inclusion $A_p(G) \subset OA_p(G)$ is a contraction with dense range, so that
$OA_p(G)$ has also an approximate identity, say $( e_\alpha )_\alpha$, bounded by one. This fact and Theorem \ref{ccthm}
enable us to define complete contractions
\[
  \iota_\alpha \!: OA_q(G) \to OA_p(G), \quad \phi \mapsto e_\alpha \phi.
\]
Define $\iota \!: OA_q(G) \to OA_p(G)^{\ast\ast}$ as the pointwise $w^\ast$-limit of $( \iota_\alpha )_\alpha$. It follows that $\iota$ is a complete contraction.
Since there is a complete contraction from $OA_p(G)$ into ${\cal C}_0(G)$, the approximate identity $( e_\alpha )_\alpha$ for $A_p(G)$ is also a bounded approximate identity for
${\cal C}_0(G)$. It follows that $\iota_\alpha(\phi) \to \phi$ uniformly in ${\cal C}_0(G)$. Since uniform convergence in ${\cal C}_0(G)$ implies weak convergence in $OA_p(G)$, it follows
that $\iota(OA_q(G)) \subset OA_p(G)$.
\end{proof}
\par
We now give an operator version of the Leptin--Herz theorem (\cite[Theorem 10.4]{Pie}), which characterizes the amenable, locally compact groups through the existence of bounded
approximte identities in their Fig\`a-Talamanca--Herz algebras:
\begin{theorem} \label{Horst&Carl}
For a locally compact group $G$, the following are equivalent:
\begin{items}
\item $G$ is amenable.
\item For each $p \in (1,\infty)$, the Banach algebra $OA_p(G)$ has an approximate identity bounded by one.
\item For each $p \in (1,\infty)$, the Banach algebra $OA_p(G)$ has a bounded approximate identity.
\item There is $p \in (1,\infty)$ such that $OA_p(G)$ has a bounded approximate identity.
\end{items}
\end{theorem}
\begin{proof}
(i) $\Longrightarrow$ (ii): see the proof of Corollary \ref{inclcor}.
\par
(ii) $\Longrightarrow$ (iii) and (iii) $\Longrightarrow$ (iv) are trivial.
\par
(iv) $\Longrightarrow$ (i): Let $q \in (1,\infty)$ be such that $\frac{1}{p} +\frac{1}{q} = 1$. Mimicking the argument (b) in the proof of \cite[Theorem 10.4]{Pie}, we see that 
$\| \lambda_q(f) \|_{\mathrm{cb}} = \| f \|_1$ for each positive $f \in {\cal C}_0(G)$ with compact support.
By Proposition \ref{cbprop}, we have $\| \lambda_q(f) \|_{\mathrm{cb}} = \| \lambda_q(f) \|$ for each such $f$ and hence, by continuity,  $\| \lambda_q(f) \| = \| f \|_1$ for all positive $f \in L^1(G)$. 
By \cite[Theorem 9.6]{Pie}, this entails the amenability of $G$.
\end{proof}
\par
For any locally compact group $G$, the Fourier algebra $A(G) = \VN(G)_\ast$ carries a canonical operator space structure. As mentioned earlier, we have $OA(G) := OA_2(G) \cong A(G)$ as Banach spaces.
We conclude this section with an example that quenches the hope that $A(G) \cong OA(G)$ even as operator spaces:
\begin{example}
Let $G$ be a compact, non-abelian group. Then $\VN(G)$, for some $n \geq 2$, contains ${\mathbb M}_n$ as a central direct summand. Assume towards a contradiction that $A(G) \cong OA(G)$ as operator spaces.
By Proposition \ref{checkprop}, (\ref{checkmap}) is a completely isometric isomorphism of $OA(G)$. Hence, 
\[
  ^\vee \!: A(G) \to A(G), \quad \phi \mapsto \check{\phi}
\]
is also a completely isometric isomorphism. By \cite[Proposition 3.2.2]{ER}, the adjoint of $^\vee$ --- which me denote likewise by $^\vee$ --- is then a complete isometry on $\VN(G)$. 
A simple calculation, however, shows that $\check{T}$ is just the Banach space adjoint of $T$ for any $T \in \VN(G)$. In particular, if $T \in {\mathbb M}_n$, then $\check{T}$ is the transpose of $T$. 
Taking the transpose on ${\mathbb M}_n$, however, is not a complete isometry (\cite[Proposition 2.2.7]{ER}). If $\VN(G)$ contains ${\mathbb M}_n$ as a central direct summand for arbitrarily large $n \in \posints$
--- e.g.\ if $G = SO(3)$ ---, then $^\vee$ is not even completely bounded. 
\end{example}
\begin{remarks}
\item If $G$ is abelian, then $A(G)$ is the predual of the commutative $\cstar$-algebra $\VN(G)$. Hence, the canonical operator space structure on $A(G)$ is $\max A(G)$ (\cite[(3.3.13) and Proposition 3.3.1]{ER}), and
the identity map from $A(G)$ to $OA(G)$ is a complete contraction. 
\item We conjecture that the identity from $A(G)$ to $OA(G)$ is a complete contraction at least for each amenable $G$. The problem when trying to mimic the proof of Corollary \ref{inclcor} in this context
is that $L^2(G)_c$ is, in general, not an operator $2$-space: If it were an operator $2$-space then \cite[Theorems 6.5 and 6.9]{Pis2} would imply that the identity on $L^2(G)_c$ factors through $OL^2(G)$;
since, however, 
\[
  {\cal CB}(L^2(G)_c,OL^2(G)) = {\cal CB}(OL^2(G),L^2(G)_c) = {\cal S}^4_{\dim L^2(G)}
\]
by \cite[Satz 5.2.28]{Lam}, this would, in turn, imply that $\id_{L^2(G)} \in {\cal S}^4_{\dim L^2(G)}$, which is impossible for infinite $G$ (in \cite{Lam} only separable Hilbert spaces are treated, but Lambert's arguments 
work for arbitrary Hilbert spaces). Nevertheless, it is sufficient for the proof of our conjecture that the amplification
\[
  \lambda_2(M(G)) \ni \lambda_2(\mu) \mapsto \lambda_2(\mu) \tensor \id_{L^2(G)}
\]
is a complete isometry into ${\cal CB}(OL^2(G,L^2(G)_c)$: an apparently much weaker statement than $L^2(G)_c$ being an operator $2$-space. {\it If\/} our conjecture is true, than it is easy to see that $G$ is amenable if and 
only if $OA_p(G)$ is operator amenable for one (and equivalently for all) $p \in (1,\infty)$.
\end{remarks}
\section*{Conclusion}
Our operator Fig\`a-Talamanca--Herz algebras $OA_p(G)$ can be considered the appropriate operator space substitutes for the classical Fig\`a-Talamanca--Herz algebras $A_p(G)$ in the following sense: For many of the classical 
theorems on the algebras $A_p(G)$, there is an operator space counterpart for the completely contractive algebras $OA_p(G)$. Nevertheless, the definition of $OA_p(G)$ has two drawbacks:
\begin{enumerate}
\item In general, we don't have $A_p(G) \cong OA_p(G)$ as Banach spaces (even if we are willing to put up with merely topological and not isometric isomorphism).
\item The isometric isomorphism $A(G) \cong OA(G)$ holds only at the Banach space level and is, in general, no complete isomorphism.
\end{enumerate}
\par
This leaves the question open of whether there is a canonical way of assigning to each $p \in (1,\infty)$ an operator space structure on $A_p(G)$, i.e.\ with $A_p(G)$ as underlying Banach space, 
such that 
\begin{itemize}
\item $A_p(G)$ is a completely contractive Banach algebra and, 
\item for $p = 2$, we obtain $A(G)$ with its canonical operator space structure.
\end{itemize}
Since the canonical operator space structure on $A(G)$ stems from the column space structure on $L^2(G)$, this might require to extend the notion of column space to arbitrary $L^p$-spaces. 
\dated
\vfill
\renewcommand{\baselinestretch}{1.2}
\begin{tabbing} 
{\it Address\/}: \= Department of Mathematical and Statistical Sciences \\
                 \> University of Alberta \\
                 \> Edmonton, Alberta \\
                 \> Canada, T6G 2G1 \\ \medskip
{\it E-mail\/}:  \> {\tt vrunde@ualberta.ca} \\[\medskipamount]
{\it URL\/}: \> {\tt http://www.math.ualberta.ca/$^\sim$runde/}  
\end{tabbing}

\end{document}